\documentclass[12pt]{article}

\usepackage[utf8]{inputenc}

\usepackage{geometry}
\usepackage{color}
\usepackage{amsmath}
\usepackage{amssymb}
\usepackage{amsthm}

\usepackage{hyperref}

\newtheorem{Lemma}{Lemma}
\newtheorem{Theorem}{Theorem}
\newtheorem{Proposition}{Proposition}

\newtheorem{Remark}{Remark}
\newtheorem{Corollary}{Corollary}
\newtheorem{Example}{Example}

\newlength{\wdth}

\def\E{\hskip.15ex\mathrm{E}\hskip.10ex}
\def\P{\mathrm{P}}

\def\phi{\varphi}

\title{On convergence rate bounds for a class of nonlinear Markov chains}

\author{A.A. Shchegolev\footnote{National Research University Higher School of Economics, Moscow, Russian Federation; email: ashchegolev@hse.ru} \, and A.Yu. Veretennikov\footnote{Institute for Information Transmission Problems, Moscow, Russian Federation; email: ayv@iitp.ru}}

\begin{document}

\maketitle

\begin{abstract}
A new approach is developed for evaluating the convergence rate for nonlinear Markov chains (MC) based on the recently developed spectral radius technique of markovian coupling for linear MC and the idea of small nonlinear perturbations of linear MC. The method further enhances recent advances in the problem of convergence for such models. The new convergence rate may be used, in particular, for the justification of $D$-condition in the Extreme Values theory.

\medskip


\medskip


\end{abstract}

\section{Introduction}

Ergodic properties of nonlinear Markov chains (nMC)  were studied recently in  \cite{Butkovsky, Butkovsky2, Butkovsky3, Shchegolev, Shchegolev2}  
among others.  
In earlier papers on {\em linear} Markov chains, ergodic properties were studied by lots of authors; we only mention a few of the most famous and important, among which are  Markov, Kolmogorov, Doeblin, Doob, Dobrushin; see a more complete list of references in \cite{VV22}. For nonlinear MC
 -- which is another name for ``processes with complete connections'' -- see the monographs  \cite{Iosifescu}  and \cite{Kolokoltsov} with references including historic ones. The 
ergodicity assumptions given in \cite{Iosifescu} are not very explicit, and conditions from \cite{Butkovsky} and \cite{Shchegolev} may be regarded as efficient sufficient conditions.

In this paper new enhanced conditions sufficient for exponential ergodicity are offered.  The whole setting is applicable both to general and finite or countable state spaces. The new conditions are not only a bit less restrictive, but {\em may} also provide better convergence rates in comparison to those in \cite{Butkovsky} and \cite{Shchegolev}. 
We treat general state spaces, with a particular emphasis on discrete ones; see examples in section \ref{sec4} and for the ``linear'' cases in~\cite{VV22}. 

So, let $\left(X_n^\mu\right)_{n\in\mathbb{Z}_{+}}$ be a nonlinear Markov chain with a 
state space $(E, \mathcal{E})$, 
initial distribution $\text{Law}\left(X_0^\mu\right) = \mu$, $\mu \in \mathcal{P}(E)$, and transition probabilities
$\mathbb{P}_{\mu_n}(x, B) = \mathbb{P}\left(X_{n+1}^\mu \in B | X_n^\mu = x; {\cal L}(X^\mu_n) = \mu_n\right)$,
where $x \in E$, $B \in \mathcal{E}$, $n \in \mathbb{Z}_{+}$ and $\mu_n := \text{Law}(X_n^\mu)$. The calculus uses some ideas from \cite{Butkovsky, ButkovskyVeretennikov,  Shchegolev, VV20,VV22} with appropriate changes. 

Recall that in  \cite{Butkovsky} it is demonstrated that the ``usual'' Markov -- Dobrushin condition taken from the ``linear'' MC theory,
\begin{equation}\label{NonlinearDobrushinAnalog}
\sup_{\mu,\nu\in\mathcal{P}(E)}\|\mathbb{P}_\mu(x,\cdot) - \mathbb{P}_\nu(x',\cdot)\|_{TV} \le 2(1-\alpha),\quad \alpha > 0, \quad x,x'\in E,
\end{equation}
{\em does not guarantee} convergence properties. Here $\displaystyle \|\mu - \nu\|_{TV}= 2\sup_A (\mu(A) - \nu(A)) = 2 - 2\int \left(\frac{\mu(dx)}{(\mu+\nu)(dx)}\wedge\frac{\nu(dx)}{(\mu+\nu)(dx)}\right) (\mu+\nu)(dx)$ is the total variation metric between  measures $\mu$ and $\nu$, and the value of Markov -- Dobrushin's coefficient $\alpha$ for the ``nonlinear'' case is defined as 
\begin{equation}\label{adef}
\alpha:= \inf_{x,x', \mu, \nu} 
\int \left(\frac{\mathbb{P}_\mu^{}(x,dy)}{\mathbb{P}_\nu^{}(x',dy)} \wedge 1 \right) \mathbb{P}_\nu^{}(x',dy), 
\end{equation}
or, in the finite or countable state space 
\begin{equation}\label{adefin}
\alpha:=  \inf_{i,i',\mu, \nu} \sum_{j} p_{ij}(\mu)\wedge p_{i'j}(\nu), 
\end{equation}
In \cite{Butkovsky} it was proposed to supplement the condition (\ref{NonlinearDobrushinAnalog}) by the following additional one,
\begin{align}\label{ButkovskyLambda}
\|\mathbb{P}_\mu(x,\cdot) - \mathbb{P}_\nu(x,\cdot)\|_{TV} \le \lambda \|\mu - \nu\|_{TV},\quad \lambda\in [0,\alpha],\quad x\in E,\quad \mu,\nu \in \mathcal{P}(E).
\end{align}
One of the results in \cite{Butkovsky} then reads,
\begin{equation}\label{munnun}
\sup_{\mu_0,\nu_0\in\mathcal{P}(E)}\|\mu_n - \nu_n \|_{TV} \le 2(1-\alpha +\lambda)^n,\quad x,x'\in E, \quad n\ge 1.
\end{equation}
The condition (\ref{ButkovskyLambda}) may be called the Lipschitz condition of the transition kernel in the total variation metric with respect to the measure variable. 
Under the combination of (\ref{NonlinearDobrushinAnalog}) and (\ref{ButkovskyLambda}) the property of the uniform exponential ergodicity and existence and uniqueness of invariant measure hold in the case of $\lambda <\alpha$; in the special case of  $\lambda =\alpha$ a bound for convergence of the order $1/n$ was shown; for the case $\lambda >\alpha$ counterexamples for convergence and even for the existence of the invariant measure have been provided, see \cite{Butkovsky, Butkovsky2, Butkovsky3}. 

In \cite{Shchegolev} it was further established that the assumption $\lambda_2 <\alpha_2$ along with $\lambda\equiv \lambda_1<\infty$ also suffices for the exponential convergence, where $\alpha_2$ and $\lambda_2$ correspond to the two-step transition kernel $\mathbb{P}^{(2)}(x,dy)$ instead of $\mathbb{P}(x,dy)$ in (\ref{NonlinearDobrushinAnalog}) and (\ref{adef}) respectively, that is, 
$$
\alpha_2:= \inf_{x,x'} \int \mathbb{P}_\mu^{(2)}(x,dy)\wedge \mathbb{P}^{(2)}_\nu(x',dy) 
\equiv \inf_{x,x',\mu,\nu} \int \left(\frac{\mathbb{P}_\mu^{(2)}(x,dy)}{\mathbb{P}_\nu^{(2)}(x',dy)} \wedge 1 \right) \mathbb{P}_\nu^{(2)}(x',dy),
$$
accompanied by the assumption 
$$
\|\mathbb{P}^{(2)}_\mu(x,\cdot) - \mathbb{P}^{(2)}_\nu(x,\cdot)\|_{TV} \le \lambda_2 \|\mu - \nu\|_{TV}, \quad  x\in E,\quad \mu,\nu \in \mathcal{P}(E).
$$
We highlight that the kernel $\mathbb{P}^{(2)}_\mu(x,\cdot)$ is given by the expression (complete probability formula)
$$
\mathbb{P}^{(2)}_\mu(x,\cdot) 
= \delta_x \mathbb P_{\mu} \mathbb P_{\mu_1}, \quad \text{where} \quad \mu_1 = \mu \mathbb P_\mu,
$$
and not by the the application of the kernel $\mathbb P_{\mu}(x,\cdot)$ with a fixed $\mu$ twice; this comment also relates to the notation $\mathbb{P}_\mu^{(k)}(x,dy)$ in what follows. 
Further, according to \cite{Shchegolev2}, a similar condition $\lambda_k <\alpha_k$ for a $k$-step kernel suffices for the exponential convergence with any fixed $k$, where $\alpha_k$ is the analogue of $\alpha\equiv \alpha_1$ for the $k$-step transition kernel, that is, 
\begin{equation}\label{alphak}
\alpha_k:= 
\inf_{x,x',\mu,\nu} \int \left(\frac{\mathbb{P}_\mu^{(k)}(x,dy)}{\mathbb{P}_\nu^{(k)}(x',dy)} \wedge 1 \right) \mathbb{P}_\nu^{(k)}(x',dy), 
\end{equation}
where $\mathbb{P}_\mu^{(0)}(x,dy) = \delta_x(dy)$, $\mu_0=\mu$,  and by induction
\begin{equation}\label{Pmuk}
\mathbb{P}_\mu^{(k)}(x,dy) \equiv \delta_x \mathbb{P}_\mu^{(k)} (dy) := \delta_x \mathbb P_\mu^{(k-1)} \mathbb P_{\mu_{k-1}}(dy), 
\quad \mu_k := \mu_{k-1}\mathbb P_{\mu_{k-1}}, \quad k\ge 1;
\end{equation}
or, in the case of a finite state space $S$,
$$
\alpha_k:= \inf_{i,i',\mu, \nu} \sum_{j} p^{(k)}_{ij}(\mu)\wedge p^{(k)}_{i'j}(\nu), 
$$
where transition probabilities $p^{(k)}_{ij}(\mu)$ are constructed in the same way as general kernels in (\ref{Pmuk}),
and $\lambda_k$ is the analogue of $\lambda\equiv \lambda_1$ for the $k$-step transition kernel, that is, the value $\lambda_k$ in Lipschitz constant in the assumption for the $k$-step transition kernel with respect to the measure:
\begin{align}\label{Lambda_k}
\|\mathbb{P}^{(k)}_\mu(x,\cdot) - \mathbb{P}^{(k)}_\nu(x,\cdot)\|_{TV} \le \lambda_k \|\mu - \nu\|_{TV}, \quad x\in E,\quad \mu,\nu \in \mathcal{P}(E), 
\end{align}
with the complementary bound
\begin{align}\label{Lambda_k2}
\lambda_k\in [0,\alpha_k), \quad \text{with} \quad \max_{i\le k} \lambda_i <\infty.
\end{align}
Under the assumption (\ref{Lambda_k})--(\ref{Lambda_k2}) in addition to 
\begin{equation}\label{alpha_k_Shch}
\sup_{\mu,\nu\in\mathcal{P}(E)}\|\mathbb{P}^{(k)}_\mu(x,\cdot) - \mathbb{P}^{(k)}_\nu(x',\cdot)\|_{TV} \le 2(1-\alpha_k), \quad  x, x'\in E,\quad \mu,\nu \in \mathcal{P}(E),
\end{equation}
in \cite{Shchegolev2} the following exponential convergence bound has been established with some $C_k<\infty$, 
\begin{equation*}
\sup_{\mu_0,\nu_0\in\mathcal{P}(E)}\|\mu_n - \nu_n \|_{TV} \le C_k(1-\alpha_k +\lambda_k)^{[n/k]},\quad x,x'\in E, \quad n\ge 1.
\end{equation*}
along with examples that the condition $\lambda_k<\alpha_k$ may take place (for $k=2,3$) while $\lambda>\alpha$. This evidently enhances the condition of Butkovsky's theorem. Condition (\ref{alpha_k_Shch}) with any $k$ will be called MD condition, for Markov -- Dobrushin.

Let us mention that the equality $\lambda_k = \alpha_k$ also leads to some convergence no slower than with the rate $1/n$ (see \cite{Shchegolev2}); however, we will not discuss this case here.

In this paper a new idea is exploited, which replaces conditions (\ref{NonlinearDobrushinAnalog}) or (\ref{Lambda_k}) based on Markov--Dobrushin's ergodic coefficients $\alpha$, or $\alpha_k$; this idea is not explicitly related to the characteristic $\alpha$ or, respectively, $\alpha_k$, but is  based on the spectral radius approach from the recent papers on linear markovian  models \cite{ButkovskyVeretennikov, VV20, VV22}. The realisation of this idea for nMC does use some calculus and results from  \cite{Shchegolev2}; in this respect, the present paper may be considered as a further enhancement of \cite{Shchegolev2}. At the same time, the result of the corollary \ref{Cor1}  following the main theorem links the asymptotics of the convergence of the model obtained as a small perturbation of the underlying linear one to the spectral characteristic $r(\bar V)$ of this linear model (see (\ref{newV})). 
In fact, the new condition implies the MD condition for the transition kernel  $\mathbb P^{(k)}$  {\em for some $k$}; the point is that the value of this $k$ is not known in advance.

Let us also mention that in the {\em Extreme Values} (EV) theory an important role belongs to the so-called $D$-condition: under this condition, the extreme index of a stationary sequence of random variables exists, which  is a highly important feature. This condition (see \cite[Section II.3.2]{Leadbetter}) automatically follows if the sequence under consideration is exponentially mixing uniformly with respect to the initial distribution, which is exactly what the main result of this paper claims. Further links between exponentially fast ergodicity and the {\em extremal index} (see \cite[Section 2.3.2]{Leadbetter}) of a stationary (linear) Markov sequence were studied, in particular, in \cite{Roberts}, especially in relation to certain MCMC algorithms. Although this goal is not pursued in this paper directly, the results of this paper apparently indicate that there may exist similar relations of the ergodicity of nMC  with the EV theory, too. 
This is the motivation for the present article to be in this issue of the Journal. One still open question which may be also of interest to the theory of Extreme Value will be mentioned after the statements of the main results.

Note that under some other conditions the established bounds for  convergence could be much weaker (with the rate $1/n$); however, they remain uniform, which means that the $D$-condition is still applicable. Let us add that the point is the possibility of the equality 
$
\alpha_i = 0, \; \text{or,}\; \alpha_i \approx 0,
\; \forall \, 1\le i \le k, 
$
for any fixed $k$, while $r(\bar V)<1$ (see what follows): see examples in the last section, as well as the examples in \cite{VV22} (the latter only relate to the linear MC case).

The paper consists of four sections: this Introduction, Assumptions and lemmata, Main result with its proof, and Examples. 

\section{Assumptions and lemmata}\label{sec2}
Two main assumptions will be made. The {\bf first one} is that there exist a homogeneous transition kernel $\bar {\mathbb{P}}_{x}$ and a constant  $\gamma>0$ such that for each $\mu,x$ the kernels 
$\bar {\mathbb{P}}_{x}$ are absolutely continuous with respect to $\mathbb{P}_{\mu,x}$ and vice versa, so that 
$\bar {\mathbb{P}}_{x}(1,dy) \sim {\mathbb{P}}_{\mu,x}(dy)$,  
and, moreover, 
\begin{equation}\label{ppt}
\frac{\bar {\mathbb{P}}_{x}(1,dy)}{{\mathbb{P}}_{\mu,x}(dy)}  \le (1+\gamma), \quad \forall \, x,y,\mu,
\end{equation}
In what follows, the value $\gamma$ will be assumed small enough. 

The {\bf second assumption} is Lipschitz condition (\ref{Lambda_k})
for the kernel $\P_\mu(x,\cdot)$ with respect to the measure $\mu$, yet, without the complementary (\ref{Lambda_k2}). 

\begin{Lemma}\label{Lemlambdak}
The inequality (\ref{ButkovskyLambda}) implies\footnote{The question whether or not it is possible that $\lambda_k<\infty$ while $\lambda_1=\infty$ remains, in general, open.} that for any $k$, $\lambda_k<\infty$, and, moreover\footnote{More generally, for any $\lambda_1$ the bound reads as $\lambda_k \le C_k  (\lambda_1\vee \lambda_1^k)$; the most essential for the results is that for any $k$, $\lambda_k = o(1)$ as $\lambda_1\to 0$.} (for what follows assume $\lambda\equiv \lambda_1\le 1$),
\begin{equation}\label{lambda2^k-1}
\lambda_k \le C_k  \lambda_1, \quad k\ge 1, 
\quad C_k = \frac53 \, 4^k - \frac23.
\end{equation}

\end{Lemma}

\noindent
{\em Proof.}\\
{\bf 1.} Indeed, for $k=1$ we have $C_1=1$; this is the induction base. 
Also, note that by the triangle inequality
\begin{align}\label{used}
\|\mu_1 - \nu_1\| \equiv
\|\mu P_\mu - \nu P_\nu\| \le \|\mu P_\mu - \nu P_\mu\| 
+ \|\nu P_\mu - \nu P_\nu\| 
 \nonumber \\ \\ \nonumber
\le  \|\mu  - \nu\| + \|P_\mu - P_\nu\| \le (1+\lambda)\|\mu  - \nu\| 
\le 2(1\vee \lambda)\|\mu  - \nu\|. 
\end{align}

\noindent
{\bf 2.} Now, to justify the induction step, let us assume that 
\begin{align*}
\|P_\mu^{(k)}(x, \cdot) - P_\nu^{(k)}(x, \cdot)\| 
 \le C_k \lambda \|\mu  - \nu\|.
\end{align*}
Consider 
\begin{align*}
|P_\mu^{(k+1)}(x, A) - P_\nu^{(k+1)}(x, A)| 
=  | \iint 1(y\in A) (P^{}_\mu(x,dx')\underbrace{P^{(k)}_{\mu P_{\mu}}(x',dy)}_{=P^{(k)}_{\mu_1}(x',dy)} -P_\nu(x,dx')\underbrace{P^{(k)}_{\nu P_{\nu}}(x',dy)}_{=P^{(k)}_{\nu_1}(x',dy)}|
 \\\\
\le \|P^{(k)}_{\mu_1} - P^{(k)}_{\nu_1}\|
+ \|P^{}_{\mu} - P^{}_{\nu}\| 
\le  C_k \lambda \underbrace{\|\mu_1  - \nu_1\|}_{\le 2(1 \vee \lambda) \|\mu  - \nu\|} + \lambda \|\mu  - \nu\|
\le (2C_k + 1) \lambda \|\mu  - \nu\|.
\end{align*}
So, 
\begin{align*}
\|P_\mu^{(k+1)}(x, \cdot) - P_\nu^{(k+1)}(x, \cdot)\| 
 \le 2(2C_k + 1) \lambda \|\mu  - \nu\|.
\end{align*}
Hence, 
\begin{equation}\label{Ckk+1}
C_{k+1}\le (4C_k + 2).  
\end{equation}
Since $C_1 = 1$, 
we get by induction for the constant in (\ref{lambda2^k-1})
\begin{equation}\label{Ck}
C_k \le \frac53 \, 4^k - \frac23.
\end{equation}
Induction step verification: assuming (\ref{Ck}) holds for $k$, using (\ref{Ckk+1}) we compute
\begin{align*}
C_{k+1} \le 4C_k + 2 = 4\times  (\frac53 \, 4^k - \frac23) +2
= \frac53 \, 4^{k+1} - \frac83 + 2 = \frac53 \, 4^{k+1} - \frac23,
\end{align*}
as required. Hence, indeed, by induction (\ref{Ck}) holds true for the constants $C_k$ in (\ref{lambda2^k-1}). Lemma \ref{Lemlambdak} is proved. \hfill QED

~

To state the main result, it is necessary to recall the markovian coupling construction for the linear Markov chains from \cite{ButkovskyVeretennikov, VV20, VV22} in order to introduce the operator $\bar V$, which spectral radius will play a crucial role in what follows.

Let us consider transition densities with respect to some dominating measure\footnote{As it is explained in \cite{VV22}, the existence of the common dominating measure is not necessary for applying this coupling algorithm; it suffices to use such a measure for each particular quadruple $(x^1, x^2, x^3, x^4)$ which always exists. We use a unique measure $\Lambda$ in what follows just to simplify the presentation.} $\Lambda(dy) = \Lambda(dy^1 dy^2 dy^3 dy^4)$ for a coupling construction related to the transition kernel $\bar {\mathbb{P}}(x,dy)$; we highlight that they do not depend on marginal measures:
\begin{equation}\label{process_etabar}
\bar\phi(x,y):=\bar\phi_1(x,y^1)\bar\phi_2(x,y^2)\bar \phi_3(x,y^3) \bar\phi_4(x,y^4), 
\end{equation}
and correspond to the homogeneous ``unperturbed'' Markov process $\bar X = (\bar X^1, \bar X^2)$ where each component is a Markov chain with the same transition kernel $\bar {\mathbb{P}}$.  
If  $0<\kappa(x^1,x^2)<1$, then for $x = (x^1,x^2,x^3,x^4)$. These densities are with respect to an appropriate dominated measure (see \cite{VV22}):
\begin{align}\label{coupldens}
\!\!\displaystyle \bar \phi_{1}(x,u)\!:=\!\frac{\bar p_{}(x^1,u)\!-\!\bar p_{}(x^1,u) \!\wedge \!\bar p_{}(x^2,u)}{1\!-\!\bar\kappa_{}(x^1,x^2)}1(x^4\!=\!1) 
+ \bar p_{}(x^3,u)\,1(x^4\!=\!0), 
 \\\nonumber\\
\!\!\displaystyle \bar \phi_{2}(x,u)\!:=\!\frac{\bar p_{}(x^2,u)\!-\!\bar p_{}(x^1,u)\!\wedge \!\bar p_{}(x^2,u)}{1\!-\!\bar \kappa_{}(x^1,x^2)}
1(x^4\!=\!1) 
+ \bar p_{}(x^3,u)\,1(x^4\!=\!0), 
 \label{phi_12}
 \\\nonumber\\
\displaystyle \bar \phi_{3}(x,u):=
1(x^4=1)\frac{\bar p(x^1,u)\wedge
 \bar p(x^2,u)}{\bar \kappa_{}(x^1,x^2)}
 + 1(x^4=0)\bar p(x^3,u),
  \label{phi_3}
  \\\nonumber\\ 
\displaystyle \bar \phi_{4}(x,u):=
\label{phi_400} 
1(x^4=1)\left(\delta_1(u)(1-\bar\kappa_{}(x^1,x^2))+ 
\delta_0(u)\bar\kappa_{}(x^1,x^2)\right) 
+1(x^4=0)\delta_0(u),
\end{align}
where
\begin{align*}
\bar \kappa(x^1,x^2) = 
\int \left(\frac{\bar {\mathbb{P}}_{x'}(1,dy)}{\bar {\mathbb{P}}_{x}(1,dy)}\wedge 1 \right)\bar {\mathbb{P}}_{x}(1,dy),
\end{align*}
The equality $x^4=0$ signifies coupling for $\bar X^1, \bar X^2$ already realised at the previous step(s) (and, hence,  $x^4=0$ implies $x^1=x^2$ for the linear Markov chain $(\bar X^1, \bar X^2)$), while $u=0$ means a successful coupling at the present step. Note that in examples the transition matrices {\em may} be easily constructed in such a way that the function $\kappa(x,x')$ is always positive for any couple $x'\neq x$. However, for generality, in the degenerate cases these densities could be set up as follows: 
where $\bar \kappa_{0}=0$ (impossible coupling at the zero step), then let 
$$
\bar \phi_{3}(x,u):= \bar p(x^1,u) 
$$ 
instead of (\ref{phi_3}); 
and if  $\bar \kappa_{0}=1$, then let  
\begin{equation}\label{after}
\bar \phi_{1}(x,u)=\bar p(x^1,u), \quad 
\bar \phi_{2}(x,u):= \bar p(x^2,u). 
\end{equation}
instead of (\ref{coupldens}) and (\ref{phi_12}), respectively. The formula (\ref{phi_400}), which determines \(\bar \phi_{4}(x,u)\) can be accepted in all cases. 

\medskip

Random variables $\left(\eta^1_0,\eta^2_0,\xi_0,\zeta_0\right)$ are chosen according to the following joint density with respect to some dominated measure (see \cite{VV22}):
\begin{equation}\label{process_eta0}
\bar \phi_{0}(y):=\bar \phi_{1}(y^1)\bar \phi_{2}(y^2)\bar \phi_{3}(y^3) \bar \phi_{4}(y^4),
\end{equation}
for $y = (y^1,y^2,y^3, y^4)$, where, in the case of 
$0<\bar \kappa_{0}<1$, 
\begin{align}
\!\!&\displaystyle \bar \phi_{1}(u)\!:=\!\frac{p^1_{}(u)\!-\!p^1_{}(u) \!\wedge \!p^2_{}(u)}{1\!-\!\bar \kappa_{0}},
\label{phi01}
 \\\nonumber\\
\!\!&\displaystyle \bar \phi_{2}(u)\!:=\!\frac{p^1_{}(u)\!-\!p^1_{}(u)\!\wedge \!p^2_{}(u)}{1\!-\!\bar \kappa_{0}},
\label{phi_120}
 \\\nonumber\\
&\displaystyle \bar \phi_{3}(u):=\frac{p^1_{}(u)\wedge
 p^2_{}(u)}{\bar \kappa_{0}},
\label{phi_30}
 \\\nonumber\\ 
&\displaystyle \bar  \phi_{4}(u):=\left(\delta_1(u)(1-\kappa_{0})+ 
\delta_0(u)\bar \kappa_{0}\right) 
\label{phi_4000}.
\end{align}
Here
$$
\bar \kappa_0 =  
\int \left(\frac{\mu_0(dy)}{\nu_0(dy)}\wedge 1 \right)\nu_0(dy).
$$
In the degenerate cases where $\bar \kappa_{0}=0$ (impossible coupling at the zero step), then let 
$$
\bar \phi_{3}(u):= p^{1}(u) 
$$ 
instead of (\ref{phi_30}); 
and if  $\bar \kappa_{0}=1$, then let  
\begin{equation}\label{after2}
\bar \phi_{1}(u)=\bar \phi_{2}(u):= p^1(u). 
\end{equation}
instead of (\ref{phi_120}). The formula (\ref{phi_4000}), which determines \(\bar \phi_{4}(u)\) can be accepted in all cases. 
Further, the homogeneous markovian coupling algorithm is defined by the formula
\begin{align}\label{coupling}
\widetilde X^1_n:=\eta^1_n 1(\zeta_n=1)+\xi_n 1(\zeta_n=0), \quad 
\widetilde X^2_n:=\eta^2_n 1(\zeta_n=1)+\xi_n 1(\zeta_n=0).
\end{align}

\begin{Lemma}[\cite{VV22}]{\label{Lem1}}
The couple $\tilde X = (\tilde X^1, \tilde X^2)$ defined in (\ref{coupling}) along with the transition densities in (\ref{coupldens})--(\ref{after}) correspond to a Markov chain with components $\tilde X^1$ and $ \tilde X^2$, each of which is also a Markov chain equivalent to $\bar X$:
$$
\tilde X^i \sim \bar X^i, \quad i=1,2.
$$ 
\end{Lemma}

\medskip

\noindent
Consider a non-negative operator
\begin{equation}\label{newV}
\bar V^{}h(x^1,x^2):= (1 - \bar\kappa^{}(x^1,x^2)) 
\mathbb E_{x^1,x^2}h(\bar X_{1}),  
\end{equation}
where $\bar X_{n} = (\bar X^1_{n},\bar X^2_{n})$ is a couple of Markov processes 
starting from $x^1_0$ and $x^2_0$, respectively, with transition densities defined in (\ref{coupldens} -- \ref{after}) above. The notation $r(\bar  V^{})$ is used for the spectral radius of $\bar  V^{}$. Note that it coincides with the spectral radius of the projection $\hat V$ of this operator on the space $S\times S \setminus \text{diagonal}$ (see \cite{VV22}); this will not be used in this paper, but could be helpful for computing the spectral radius because of a reduction of dimension. 

\begin{Proposition}[\cite{VV22}]\label{Pro1}
Under the condition (\ref{ppt}),  for any $\delta>0$, for all $n$ large enough 
\begin{equation}\label{pro1-2}
\|\underbrace{\bar {\mathbb{P}}_{\mu^1_0}(n,\cdot)}_{=: \mu^1_n} - \underbrace{\bar {\mathbb{P}}_{\mu^2_0}(n,\cdot)}_{=: \mu^2_n}\|_{TV} 
\le 2\,\bar {\mathbb{P}}_{\mu^1_0, \mu^2_0}(\tilde X^1_n \neq \tilde X^2_n)
\le \left(r(\bar  V^{})+\delta\right)^n \|\mu^1_0 - \mu^2_0\|_{TV}.
\end{equation}
and
\begin{align}\label{newrate13c}
\limsup\limits_{n\to\infty} \frac1n \ln \|\underbrace{\bar {\mathbb{P}}_{\mu^1_0}(n,\cdot)}_{=: \mu^1_n} - \underbrace{\bar {\mathbb{P}}_{\mu^2_0}(n,\cdot)}_{=: \mu^2_n}\|_{TV}
\le \ln \left(r(\bar  V^{})\right).
\end{align}
The limit in (\ref{newrate13c}) is uniform with respect to the initial measures $\mu^1_0, \mu^2_0$.
\end{Proposition}
\noindent
{\em Proof.}
In fact, the second inequality on the following line 
$$
\frac12\,\|\underbrace{\bar {\mathbb{P}}_{\mu^1_0}(n,\cdot)}_{=: \mu^1_n} - \underbrace{\bar {\mathbb{P}}_{\mu^2_0}(n,\cdot)}_{=: \mu^2_n}\|_{TV} 
\le \bar {\mathbb{P}}_{\mu^1_0, \mu^2_0}(\tilde X^1_n \neq \tilde X^2_n)
\le\left(r(\bar  V^{})+\delta\right)^n 
$$
(while the first inequality is automatically satisfied for any coupling) for large enough $n$ has been established in \cite{VV22} without the multiplier $\|\mu^1_0 - \mu^2_0\|_{TV}$. The arguments there do allow this multiplier without any new calculus due to the ``forgotten'' term 
$1-\kappa(0)$ which was discussed in the beginning of section 2.2 of \cite{VV22}; it  {\em may} be added as a multiplier in the inequality (20) without any additional explanation, with the remark that $1-\kappa(0) = \|\mu^1_0-\mu^2_0\|/2$ (note that, in general, we cannot claim the same for $n>0$). Note that in \cite[inequality (20)]{VV22} it was stated in the form with $\mu^1_0 = \delta_{x^1}$ and $\mu^2_0 = \mu$, which does not change the reasoning and the inequality itself that could be rewritten in the following version (here we use notations $\eta^j_i$ from \cite{VV22}): 
$$
\P_{\mu^1_0,\mu^2_0}(\widetilde X^1_n\neq \widetilde X^2_n) 
\le (1-\kappa(0))\E_{\mu^1_0,\mu^2_0}\prod_{i=0}^{n-1} (1-\kappa(\tilde X^1_i,\tilde X^2_i))
\le  (1-\kappa(0))\E_{\mu^1_0,\mu^2_0}\prod_{i=0}^{n-1} (1-\kappa(\eta^1_i,\eta^2_i))
$$
with a generic initial distribution $\mu^1_0$ for the first component.
Hence, both (\ref{newrate13c}) and (\ref{pro1-2}) are proved. \hfill QED

~

\noindent
\begin{Remark}
Naturally, the estimate (\ref{newrate13c}) makes sense if 
\begin{equation}\label{sense1}
r(\bar  V^{}) < 1.
\end{equation}
For some -- actually, for many -- classes of processes the value $r(\bar  V^{})$ can be a better estimate than $1-\bar\alpha$, namely, iff 
\begin{equation}\label{sense20}
r(\bar  V^{}) < 1- \bar\alpha,
\end{equation} 
see \cite[Examples]{VV22}. 
Note that $r(\bar  V^{}) \le 1- \bar\alpha$ in all cases. 
The inequality (\ref{sense20}) is well possible in the case of small enough $\gamma$ and $\alpha \approx 0$, or even just $\alpha = 0$,  combined with a small enough value of $r(\bar V)$. Moreover, in examples it may also well occur that $r(\bar V)^2 < 1-\alpha_2$, etc., see section \ref{sec4}. 

\end{Remark}

In the next lemma it is assumed that the norm of the vector (function) on $S^2$ is a sup-norm, and $\|\bar  V\|$ is the operator norm of $\bar  V$. 
\begin{Lemma}\label{T1}
For any $\delta >0$ there exists $n_0$ such that for any $n\ge n_0$
\begin{equation}\label{sense22}
\|\bar V^{n} \| \le (r(\bar  V^{}) + \delta)^{n}.  
\end{equation} 

\end{Lemma}

\noindent
{\em Proof.} The bound (\ref{sense22}) for large enough values of $n$ follows straightforwardly from Gelfand's formula and the definition of the spectral radius and its eigenfunction. \hfill QED

\begin{Remark}\label{openq}
Note that the standing assumptions (\ref{NonlinearDobrushinAnalog}) and (\ref{ButkovskyLambda}) of the theory for nMC developed earlier in \cite{Butkovsky} and \cite{Shchegolev} -- which are weakened in some sense in this paper -- may be called compactness conditions, since the most natural way to guarantee them could be to assume that the state space (or the main auxiliary operator $\bar V$, see in what follows) is compact, among some others. One particular case of non-compact state space {\em for continuous time} counterpart of this theory, namely, for McKean -- Vlasov SDEs with an additive Wiener process was tackled in \cite{Butkovsky}. So far, the same question of convergence rates et al.  for the ``non-compact'' discrete time situations remains open. We mention it here because this issue is, of course, of interest to the EV theory.
\end{Remark}

\section{Main result}\label{sec3}

\begin{Theorem}\label{T1}
Under the assumptions  (\ref{ppt}) and (\ref{Lambda_k}) for all $k\ge 1$, 
for any $\delta >0$ 
small enough there exist $\lambda(\delta)>0$, and $\gamma(\delta)>0$, and a constant $C>0$ such that for any $\gamma<\gamma(\delta)$ and $\lambda_1<\lambda(\delta)$ and for $n$ large enough
\begin{equation}\label{rdl}
\|\mu_n - \nu_n\|_{TV} \le 2C(r(\bar V)+\delta)^n.
\end{equation}

\end{Theorem}
Some information about the constant $C$ in the bound (\ref{rdl}) 
will be available in the proof of the theorem. For some additional information in a special ``compact'' case see remark \ref{rem--after-thm} after the proof of the theorem. 

\begin{Corollary}\label{Cor1}
Under the assumptions  (\ref{ppt}) and (\ref{Lambda_k}) 
\begin{equation}\label{limln}
\limsup_{\gamma, \lambda_1\to 0}\limsup_{n} 
\frac1{n} \ln \|\mu_n - \nu_n\|_{TV} \le \ln r(\bar V).
\end{equation}
\end{Corollary}

~

\noindent
{\em Proof of theorem.}
Let $\delta>0$. First of all, in both cases it suffices to consider the situation where $(r(\bar V)  +\delta)^k < 1-\alpha_k$ for each $k$; otherwise the theorem folows from \cite[Theorem 1]{Shchegolev2}.

\medskip

Now assuming $(r(\bar V)  +\delta)^k < 1-\alpha_k$ for each $k$, we are going to verify the inequality similar to  (\ref{alpha_k_Shch}) {\bf with another constant} which replaces $1-\alpha_k$ for some finite $k$ by the value $(r(\bar V)+\delta)^k$. Then, under the assumption of (\ref{Lambda_k}), due to the inequality (\ref{lambda2^k-1}) the desired bound will follow from \cite[Theorem 1]{Shchegolev2} for any $\lambda_1$ small enough; indeed, small $\lambda_1$ implies small $\lambda_k$ for each particular value of $k$. In fact, we shall see in the end that the inequality $(r(\bar V)  +\delta)^k < 1-\alpha_k$ {\bf for each} $k$ is not possible. Yet, the point is that it is  not known in advance for which $k$  the latter inequality fails; hence, the new bound for the rate of convergence under the condition that both $\lambda$ and $\gamma$ are small enough is based on the unique characteristic $r(\bar V)$ rather than on any $1-\alpha_k$.

To achieve the inequality (\ref{alpha_k_Shch}) with a new constant instead of $1-\alpha_k$, let us estimate the distance $\|\mu_n - \nu_n\|_{TV}$ for two given initial distributions $\mu_0$ and $\nu_0$ (here $n$ will be chosen and fixed a bit later depending on $\delta$)  by using the triangle inequality applied twice,
\begin{align}\label{2triangle}
\|\mu_n - \nu_n\|_{TV} \le \|\mu_n - \bar \mu_n\|_{TV} + \|\nu_n - \bar \nu_n\|_{TV} + \|\bar \mu_n - \bar \nu_n\|_{TV}.
\end{align}
Here we assume that $\bar \mu_0 =  \mu_0$ and $\bar \nu_0 =  \nu_0$.  
According to proposition \ref{Pro1}, 
\begin{align}\label{ndelta}
\|\bar \mu_{n} - \bar \nu_{n}\|_{TV} \le 2(r(\bar V) + \delta)^{n} 
\end{align}
for any $n$ large enough, independently of the initial distributions\footnote{Even $\|\bar \mu_{n} - \bar \nu_{n}\|_{TV} \le (r(\bar V) + \delta)^{n}\|\bar \mu_{0} - \bar \nu_{0}\|_{TV} = (r(\bar V) + \delta)^{n}\|\mu_{0} - \nu_{0}\|_{TV}$, although, in the limit this does not affect the result; yet, if the distance $\|\mu_{0} - \nu_{0}\|_{TV}$ may be evaluated, this could decrease the value of $n_\delta$ in what follows.}. 

\medskip

\noindent
Now given $\delta>0$, let us  choose
\begin{equation}\label{ndeltadef}
n_\delta:= \inf(n: \text{(\ref{ndelta}) holds for this $n$}).
\end{equation}
Note that the set under the $\inf$ symbol is not empty.
Further for this fixed value of $n=n_\delta$ let us estimate the values of $\|\mu_n - \bar \mu_n\|_{TV}$ and $\|\nu_n - \bar \nu_n\|_{TV}$: show that both values are arbitrarily small if $\gamma$ is small enough. Both terms are quite similar, so we only consider the first one, $\|\mu_n - \bar \mu_n\|_{TV}$. Let
\begin{equation}\label{rng}
\rho_n : = \prod_{i=0}^{n-1} \frac{\bar p_{}(\tilde X^1_i, \tilde X^1_{i+1})}{p_{\mu_i}(\tilde X^1_i, \tilde X^1_{i+1})} 
\le (1+\gamma)^{n},
\end{equation}
where $\gamma$ in the constant from the condition (\ref{ppt}).
Here $\tilde X^1$ is the observable trajectory of our nonlinear MC. (Note that the second component $\tilde X^2$ is not used and is not needed here because we are to compare just the measures $\mu_n$ and $\bar \mu_n$ which both only relate to one component of the pair $(\tilde X^1, \tilde X^2)$.) 
Moreover, $\rho_n$ serves as a probability density of the measure ${\mathbb{P}}$ with respect to ${\mathbb{P}}^{\rho_n}$  and vice versa, $(\rho_n)^{-1}$ is a probability density of the measure  ${\mathbb{P}}^{\rho_n}$ with respect to ${\mathbb{P}}$ on ${\cal F}_n$. 
Here ${\mathbb{E}}^{\rho_n}\xi := {\mathbb{E}} \rho_n \xi$.
Under the measure ${\mathbb{P}}^{\rho_n}$ the nMC $\tilde X^1_k, \, k\le n$ has the same distribution in the space of trajectories as the linear MC $\bar X^1_k, \, k\le n$ under the initial measure ${\mathbb{P}}$.

~

As it is well known, the total variation distance may be estimated from above via the density of one measure with respect to the other. Denote by $\mu_{0,n}$ ($\bar \mu_{0,n}$) the measure on the space of trajectories  on $[0,n]$ of the process $\tilde X^1$ (respectively, of the process $\bar X^1$). We have for $n\ge 1$,
\begin{align*}
\frac12 \|\mu_n -\bar \mu_n\|_{TV} \le \frac12 \|\mu_{0,n} -\bar \mu_{0,n}\|_{TV} 
= \sup_{A\in S^n} \int_A \left(1 - \frac{\bar \mu_{0,n}(dz)}{\mu_{0,n}(dz)}\right)\mu_{0,n}(dz)
 \\\\
 =  \int \left(1 - \frac{\bar \mu_{0,n}(dz)}{\mu_{0,n}(dz)}\wedge 1\right)\mu_{0,n}(dz)
= \int \left(1 - 1\wedge \prod_{i=0}^{n-1} \frac{\bar p_{}(\tilde x^1_i, \tilde x^1_{i+1})}{p_{\mu_i}(\tilde x^1_i, \tilde x^1_{i+1})} \right) 
\prod_{i=0}^{n-1} p_{\mu_i}(\tilde x^1_i, \tilde x^1_{i+1}) \mu_0(d\tilde x_0)
 \\\\
= 1- {\mathbb{E}}(\rho_n\wedge 1) = {\mathbb{E}}(1-\rho_n\wedge 1)
=   {\mathbb{E}}(1-\rho_n)1(\rho_n<1) 
 \\\\
\stackrel{CBS}\le \sqrt{{\mathbb{E}}(1-\rho_n)^2}
= \sqrt{1-2{\mathbb{E}}\rho_n + {\mathbb{E}}\rho_n^2} 
= \sqrt{{\mathbb{E}}\rho_n^2 - 1}
\le \sqrt{(1+\gamma)^{2n}-1} \le \sqrt{3n\gamma}, 
\end{align*} 
where the last inequality holds for $\gamma>0$ small enough, if $n$ is fixed due to the limit $\lim_{\gamma\downarrow 0}((1+\gamma)^{2n}-1)/\gamma = 2n$; ``CBS'' stands for the Cauchy -- Bunyakovsky -- Schwarz inequality. 
Similarly, for $\gamma>0$ small enough
\begin{align*}
\frac12 \|\nu_n -\bar \nu_n\|_{TV} \le \sqrt{3n\gamma}, \quad n\ge 1.
\end{align*} 
Let us choose $n=n_\delta$ and
$$
3\gamma < \frac1{16 n_\delta} \, (1-2(r(\bar V)+\delta)^{n_\delta})^{2}.
$$
Then the desired bound (\ref{rdl}) holds for $n=n_\delta$, because the value $\|\nu_n-\bar\nu_n\|$ is estimated in the same way, which due to  (\ref{2triangle}) leads to the bound 
\begin{align*}
\|\mu_{n_\delta} - \nu_{n_\delta}\|_{TV} 
\le  \|\bar \mu_{n_\delta} - \bar \nu_{n_\delta}\|_{TV} 
+ 4 \sqrt{3n_\delta\gamma} 
 \\\\
\le  2(r(\bar V) + \delta)^{n_\delta} + 1-2(r(\bar V)+\delta)^{n_\delta} = 1 <2. 
\end{align*}
This serves as an analogue of Markov -- Dobrushin's constant $2(1-\alpha_{n_\delta})$ for $n_\delta$ steps, as required, with $\alpha_{n_\delta}$ replaced by  $1/2$. 
Now by virtue of \cite[Theorem 1]{Shchegolev2} and of proposition \ref{Pro1} 
we obtain the inequality (\ref{rdl}) for all values of $n$. The constant $C$ here can be chosen as 
\begin{equation}\label{C1}
C=\max_{i\le n_\delta}(1+\lambda_i).
\end{equation}
The theorem follows. \hfill QED

\medskip

\begin{Remark}\label{rem--after-thm}
Under the same assumptions of the theorem, if, in addition, the operator $\bar V$ is compact\footnote{which is always the case for finite matrices}, then there exist $\lambda(\delta)>0$ and $\gamma(\delta)>0$ and $C\ge 1$ such that for any $\gamma<\gamma(\delta)$ and $\lambda_1<\lambda(\delta)$
and for all $n\ge 1$
\begin{equation}\label{rl}
\|\mu_n - \nu_n\|_{TV} \le 2C (r(\bar V)+ \delta)^n.
\end{equation}
\noindent
Indeed, 
in the case where the operator $\bar V$ is of the Frobenius type, that is, compact and evidently non-negative, a better non-asymptotic bound may be guaranteed:
\begin{align}\label{Cn}
\|\bar \mu_n - \bar \nu_n\|_{TV} \le 2C \bar r^n.
\end{align}
It follows from the bound
$$
\|\bar \mu_n - \bar \nu_n\|_{TV} \le 2\P(\tilde X^1_n\neq \tilde X^2_n) \le 2\bar V^n 1(x), 
$$
and from the identity 
$$
\bar V e (x) = \bar r e(x), \quad \& \quad \bar V^n e (x) = \bar r^n e(x), 
$$
and from the double inequality
$$
0<\min_x e(x) \le e(x) \le \max_x e(x), \quad \forall \, x.
$$
Indeed, all of the above implies straightforwardly that 
$$
\bar V {\bf 1} (x) \le c^{-1} \bar V e (x), \quad \& \quad 
 \bar V^n {\bf 1} (x)\le c^{-1} \bar V^n e (x) \le \frac{C}{c} \bar V^n 1 (x). 
$$
Hence, (\ref{Cn}) follows with $\displaystyle C= \frac{\max_x e(x)}{\min_x e(x)}$. 

This better inequality (\ref{Cn}) in comparison to the more general estimate (\ref{ndelta}) due to proposition \ref{Pro1} leads to the same asymptotic bound (\ref{limln}) in corollary \ref{Cor1}, as $\gamma, \lambda_1\to 0$. This is why it was not included in the statement of the theorem. Yet, it might be useful in some examples, so it is presented here as a remark. 
\end{Remark}

\begin{Remark}
One may say that the established bound (\ref{rdl}) is of the same meaning as Markov -- Dobrushin's inequality (\ref{NonlinearDobrushinAnalog}), just with ``another $\alpha$'' and an additional multiplier $C$. However, the point is that $\alpha$ is a very particular Markov -- Dobrushin's constant (see (\ref{adef}), (\ref{adefin})) and $\alpha_k$ is its analogue for $k$ steps (see (\ref{alphak})), while $r(\bar V^{})$, or $r(\bar V^{})+ \delta$ may well be less than $1-\alpha$, or $(1-\alpha_k)^{1/k}$ for any $k$ fixed in advance. 

Note that  we do not claim that $r(\bar V^{}) \le (1-\alpha_k)^{1/k}$ for all $k\ge 1$, this question is currently open. However, if the opposite inequality $r(\bar V^{}) > (1-\alpha_k)^{1/k}$ occurs, it is always possible to re-arrange the ``spectral radius'' approach by defining the analogue of the operator $\bar V$ for $k$ steps: say, denote it $\bar V^{(k)}$ (see \cite{VV22}). Then it is true that 
$r(\bar V^{(k)}) \le (1-\alpha_k)$, and likely in most examples it would be $r(\bar V^{(k)}) < (1-\alpha_k)$. 

\end{Remark}
\section{Examples}\label{sec4}
The following examples were tackled using the SymPy Python library and Wolfram Mathematica 12.3. The spectral radius was calculated by taking the modulus of the greatest eigenvalue for the matrix of the coupled process. In turn, the coupled process matrix was calculated using the  formulae \eqref{coupldens}--\eqref{phi_400}.

\begin{Example}\label{Ex1}
Let us consider a discrete nMC with the state space $(E,\mathcal{E}) = \left(\{1, 2, 3, 4\}, 2^{\{1, 2, 3, 4\}}\right)$, the initial distribution $\mu_0 = (\mu^0_1, \mu^0_2, \mu^0_3, \mu^0_4)$ and transition probability matrix:
$$P_{\mu,x} = 
\begin{pmatrix}
0.4 - \kappa \mu_{1} & 0.2 & 0.2 + \kappa \mu_{1} & 0.2\\
0.3 & 0.4 & 0.2 & 0.1\\
0.2 & 0.2 & 0.4 & 0.2\\
0.2 & 0.1 & 0.2 & 0.5
\end{pmatrix},\quad
\bar P_x = 
\begin{pmatrix}
0.4 & 0.2 & 0.2 & 0.2\\
0.3 & 0.4 & 0.2 & 0.1\\
0.2 & 0.2 & 0.4 & 0.2\\
0.2 & 0.1 & 0.2 & 0.5
\end{pmatrix},
$$
where $0 \le \kappa < 0.3$ (to ensure that $\alpha$ for $P_{\mu,x}$ can be approached for the pair of states $\{2,4\}$). Then, for the matrix $P_{\mu,x}$ we have $\alpha = 0.6$ and $\lambda = \kappa$.

Consider the matrix $\bar P_x$. The estimate $1-\bar\alpha$ for this matrix equals $0.4$. The matrix $\bar V$ for the coupled process reads, 
$$
\left[\begin{array}{cccccccccccc}
\frac{1}{10} & 0 & 0 & 0 & 0 & 0 & 0 & 0 & 0 & 0 & \frac{1}{10} & 0\\
0 & \frac{1}{5} & 0 & 0 & 0 & 0 & 0 & 0 & 0 & 0 & 0 & 0\\
0 & 0 & \frac{1}{5} & 0 & 0 & \frac{1}{10} & 0 & 0 & 0 & 0 & 0 & 0\\
0 & 0 & 0 & \frac{1}{10} & 0 & \frac{1}{10} & 0 & 0 & 0 & 0 & 0 & 0\\
0 & \frac{1}{15} & \frac{1}{30} & 0 & \frac{2}{15} & \frac{1}{15} & 0 & 0 & 0 & 0 & 0 & 0\\
0 & 0 & \frac{1}{10} & 0 & 0 & \frac{3}{10} & 0 & 0 & 0 & 0 & 0 & 0\\
0 & 0 & 0 & 0 & 0 & 0 & \frac{1}{5} & 0 & 0 & 0 & 0 & 0\\
0 & 0 & 0 & 0 & 0 & 0 & \frac{1}{15} & \frac{2}{15} & 0 & \frac{1}{30} & \frac{1}{15} & 0\\
0 & 0 & 0 & 0 & 0 & \frac{1}{10} & 0 & 0 & \frac{1}{5} & 0 & 0 & 0\\
0 & 0 & 0 & 0 & 0 & 0 & 0 & 0 & 0 & \frac{1}{5} & \frac{1}{10} & 0\\
0 & 0 & 0 & 0 & 0 & 0 & 0 & 0 & 0 & \frac{1}{10} & \frac{3}{10} & 0\\
0 & 0 & 0 & 0 & 0 & 0 & 0 & 0 & 0 & 0 & \frac{1}{10} & \frac{1}{5}
\end{array}\right]
$$
The eigenvalues for the coupled process matrix are as follows:
$$
\left[ \frac{1}{4} - \frac{\sqrt{5}}{20}, \  \frac{1}{4} - \frac{\sqrt{5}}{20}, \  \frac{\sqrt{5}}{20} + \frac{1}{4}, \  \frac{\sqrt{5}}{20} + \frac{1}{4}, \  \frac{1}{10}, \  \frac{1}{10}, \  \frac{1}{5}, \  \frac{1}{5}, \  \frac{1}{5}, \  \frac{1}{5}, \  \frac{2}{15}, \  \frac{2}{15}\right]
$$
Therefore, $r = \frac{\sqrt{5}}{20} + \frac{1}{4} \approx 0.3618$, $r^2 \approx 0.1309$, $r^3 \approx 0.04736$, $r^4 \approx 0.017135$ , $r^5 \approx 0.0061996$, $1- \bar\alpha_2 = 0.15$, $1 - \bar\alpha_3 = 0.055$, $1 - \bar\alpha_4 = 0.02$, $1 - \bar\alpha_5 = 0.00725$.
We have $r < 1 - \bar\alpha = 1 - \alpha$, $r^2 < 1-\bar\alpha_2$, $r^3 < 1-\bar\alpha_3$, $r^4 < 1-\bar\alpha_4$, $r^5 < 1-\bar\alpha_5$. This shows that the new assumption based on the spectral radius approach for linear MC and on nonlinear small perturbations is weaker that each of the MD conditions, at least, for $k=1, \ldots, 5$.
\end{Example}

\begin{Example}\label{Ex2}
Let us consider another discrete nMC with a bit more involved nonlinear components. Let the nMC have the state space $(E,\mathcal{E}) = \left(\{1, 2, 3, 4, 5, 6\}, 2^{\{1, 2, 3, 4, 5, 6\}}\right)$, the initial distribution $\mu_0 = (\mu^0_1, \mu^0_2, \mu^0_3, \mu^0_4, \mu^0_5, \mu^0_6)$ and transition probability matrix:
$$P_{\mu,x} = 
\begin{pmatrix}
0.4 - \kappa\mu_{1}& 0.2 + \kappa\mu_{1} & 0.1 & 0.1 & 0.1 & 0.1\\
0.2 & 0.3 - \kappa\mu_{2}  & 0.2 + \kappa\mu_{2} & 0.1 & 0.1 & 0.1\\
0.1 & 0.2 & 0.3 - \kappa\mu_{3} & 0.2 + \kappa\mu_{3} & 0.1 & 0.1\\
0.1 & 0.1 & 0.2 & 0.3 - \kappa\mu_{4} & 0.2 + \kappa\mu_{4} & 0.1\\
0.1 & 0.1 & 0.1 & 0.2 & 0.3 - \kappa\mu_{5} & 0.2 + \kappa\mu_{5}\\
0.1 & 0.1 & 0.1 & 0.1 & 0.2 + \kappa\mu_{6} & 0.4 - \kappa\mu_{6}
\end{pmatrix},$$
where $0 \le \kappa < 0.3$.
The transition probability matrix for the corresponding linear Markov chain is
$$
\bar P_x = 
\begin{pmatrix}
0.4 & 0.2 & 0.1 & 0.1 & 0.1 & 0.1\\
0.2 & 0.3 & 0.2 & 0.1 & 0.1 & 0.1\\
0.1 & 0.2 & 0.3 & 0.2 & 0.1 & 0.1\\
0.1 & 0.1 & 0.2 & 0.3 & 0.2 & 0.1\\
0.1 & 0.1 & 0.1 & 0.2 & 0.3 & 0.2\\
0.1 & 0.1 & 0.1 & 0.1 & 0.2 & 0.4
\end{pmatrix}.
$$
Then, for the matrix $P_{\mu,x}$ we have $\alpha = 0.6$ and $\lambda = \kappa$.

Consider the matrix $\bar P_x$. The estimate $1-\bar\alpha$ for this matrix equals $0.4$.
The matrix $\bar V$ for the coupled process has the size $30\times30$ with all eigenvalues bounded away from zero.

Therefore, $r \approx 0.3732$, $r^2 \approx 0.139282$, $r^3 \approx 0.05198$, $r^4 \approx 0.019399$, $r^5 \approx 0.007239986$, $1- \bar\alpha_2 = 0.16$, $1 - \bar\alpha_3 = 0.062$, $1 - \bar\alpha_4 = 0.0236$, $1 - \bar\alpha_5 = 0.0089$.
We have $r < 1 - \bar\alpha = 1 - \alpha$, $r^2 < 1-\bar\alpha_2$, $r^3 < 1-\bar\alpha_3$, $r^4 < 1-\bar\alpha_4$, $r^5 < 1-\bar\alpha_5$. Again, the new assumption based on the spectral radius approach for linear MC and on nonlinear small perturbations is weaker that each of the MD conditions, at least, for $k=1, \ldots, 5$. 
\end{Example}

\section*{Acknowledgements}
The section \ref{sec4} with all computations was prepared by the first author  within the framework of the HSE University Basic Research Program. For both authors this study was funded by the Russian Foundation for Basic Research grant 20-01-00575a: 
theorem \ref{T1} together with corollary \ref{Cor1}  in section \ref{sec3}  were established by the first author, and  proposition \ref{Pro1}, all the lemmata, and the construction of the markovian coupling in section \ref{sec2} are due to the second author.

\end{document}